\documentclass[letterpaper,10pt,conference]{ieeeconf}

\IEEEoverridecommandlockouts%
\usepackage{amsmath}
\usepackage{mathtools}
\usepackage{amssymb}
\usepackage{cases}
\usepackage{xcolor}
\usepackage{bm}
\usepackage{cite}
\usepackage{graphicx}
\usepackage{tikz}
\usepackage{pgfplots}
\pgfplotsset{compat=1.18}
\usepgfplotslibrary{groupplots,fillbetween}
\usetikzlibrary{calc}
\newtheorem{proposition}{Proposition}
\newtheorem{theorem}{Theorem}

\newtheorem{example}{Example}
\newtheorem{remark}{Remark}
\newtheorem{problem}{Problem}

\title{Capacity-Aware Redistribution Control for Mitigation of Cascading Failures in Networks}

\author{Ujjwal Pratap and Christian A. Hans
\thanks{Ujjwal Pratap and Christian A. Hans are with the
Automation and Sensorics in Networked Systems (ASN) Group, University of
Kassel, Kassel, Germany (e-mail: \texttt{\{ujjwal.pratap,hans\}@uni-kassel.de}).}%
\thanks{This work was partially supported by the German Federal Ministry for Economic
Affairs and Energy (BMWE) under grant number 03EN1122A (IntelliEMEs).}}

\begin{document}

\maketitle
\thispagestyle{empty}
\pagestyle{empty}

\begin{abstract}
This paper studies the control of irreversible cascading failures in networks with finite node capacities.
Each node carries a nonnegative divisible quantity, referred to as load, and
fails irreversibly when its load exceeds its capacity. The network dynamics
are represented by topology-dependent row-stochastic matrices. When nodes fail,
these matrices must be updated. However, directly deleting rows and
columns corresponding to failed nodes can reduce the total load in the network. We therefore
introduce a topology update that preserves unit row sums and quantifies any
unserved load explicitly. We derive the resulting load-balance relation and a
necessary and sufficient condition for lossless redistribution after
simultaneous failures. Based on these results, we formulate a two-stage
optimization problem: the first stage minimizes unserved load, and the second stage
minimizes the highest capacity utilization among the surviving nodes.
Simulations on grid and Barabási–Albert networks show that, unlike the standard equal and residual-capacity redistribution approaches, the proposed method prevents secondary failures in all tested cases, with a tradeoff against served load on scale-free networks.
\end{abstract}

\section{Introduction}
\label{sec:introduction}

Many networked systems, such as power, transportation, communication, and service networks, move a conserved quantity, hereafter referred to as load, among nodes with finite capacities such as substations, stations, or routers. When the load at a node exceeds its capacity, the node fails permanently, and its load is either transferred to surviving nodes or remains unserved.  Redistribution may overload further nodes and trigger a cascade of
failures~\cite{motter2002cascade,crucitti2004model,dobson2007complex}.
Each failure also changes the network topology and thereby alters subsequent
load transfers.

Capacity-based models~\cite{motter2002cascade,crucitti2004model} describe how overload and changes in network flows
can lead to propagating failures.
Decentralized redistribution using neighboring information is proposed
in~\cite{AlTakrouri2013}, and optimal reassignment of failed-link loads in
multilayer networks is studied in~\cite{kumar2021minimizing}.  Related cascade models
and mitigation strategies have been studied in the complex-networks
literature~\cite{valdez2020cascading,buldyrev2010catastrophic,motter2004control}.
For local redistribution, equal sharing can be used. It requires no optimization, but ignores
the residual capacities of surviving nodes.  It also does not coordinate
transfers from simultaneously failing nodes that share a receiver, i.e., a surviving neighbor.

In power networks, the method in~\cite{bienstock2011optimal} computes
load-shedding schedules over successive cascades using linearized
power-flow equations.  For a similar application, homogeneous load reduction, optimized shedding, and a
tree heuristic approach are compared in~\cite{pahwa2013load}.  For networks with
external inflows, \cite{Savla2014} investigates distributed routing under
irreversible link failures and reductions in link capacities.  For continuous-time flow networks, throughput and overload
properties under monotone distributed routing are derived
in~\cite{como2015throughput}.
Generally, existing cascade models typically change load only at failures
and treat each failed node separately.  Moreover, the optimization-based and
robust approaches
in~\cite{bienstock2011optimal,pahwa2013load,Savla2014,como2015throughput}
rely on power-flow models or on flow networks with external inflows, and are
thus not directly applicable when a finite node load is reassigned at
failures.  It therefore remains unclear how load transfers between
failures should adapt to node removal, and under which conditions the loads of simultaneously
failing nodes can be reassigned without overloading shared receivers.

In this work, we consider a finite divisible load held at the nodes, with no
continuing external inflow.  Nominal transfers follow a topology-dependent
row-stochastic matrix, whereas loads of failed nodes are reassigned among
surviving neighbors.  When nodes fail, directly deleting the corresponding
rows and columns can reduce the total load before redistribution. We propose an alternative to this strategy. In this context, our contributions are
as follows.  First, we show that our alternative topology update preserves row
stochasticity after single or simultaneous failures and derive the
corresponding total-load balance.  Second, when the nominal surviving-node
loads have enough capacity, we give a necessary and sufficient condition
for lossless redistribution after simultaneous failures and determine the
minimum load that must remain unserved.  Third, we formulate a two-stage
linear program that first minimizes unserved load and then minimizes the highest
surviving-node capacity utilization.  We give a necessary and sufficient condition under which no
further failures occur.  Fourth, in simulations with grid and scale-free networks, we compare our approach with standard heuristics, i.e., equal and residual-capacity redistribution, after two
initial failures with a common surviving neighbor occur, and redistribution
requires accounting for shared receiver capacities.
\section{Network Model and Failure Dynamics}
\label{sec:problem}
Consider an undirected weighted graph $G=(V,E)$ with nodes
$V=\{1,\ldots,n\}$ and edges $E \subseteq V\times V$ with nonnegative weights.  An edge $(i,j)\in E$
allows load transfer between nodes $i$ and $j$.  Load is a nonnegative
divisible quantity assigned to the nodes and can be divided among neighbors.
The fixed capacity $c_i>0$ is the maximum load that node $i$ can carry.
The load and capacity vectors are denoted by
$\bm{x}(k)=[x_1(k) ~\cdots~x_n(k)]^{\mathsf T}\in\mathbb{R}_{\geq 0}^n$ and
$\bm{c}=[c_1~\cdots~c_n]^{\mathsf T}\in\mathbb{R}_{>0}^n$,
respectively.  Let $\mathcal A(k)$ be the set of nodes active before the detection of failures at time $k$, with $\mathcal A(0)=V$.  The newly
failed and previously failed set of nodes respectively belong to
\begin{subequations}
\begin{align}
 \mathcal F(k)&=\{i\in\mathcal A(k):x_i(k)>c_i\},
 \label{eq:failure-set}\\
 \mathcal R(k)&=V\setminus\mathcal A(k).
 \label{eq:previous-set}
\end{align}
\end{subequations}

The active graph $G(k)=G[\mathcal A(k)]$ contains the active nodes and all
edges of $G$ between them.  Failures are irreversible, so
$\mathcal A(k+1)=\mathcal A(k)\setminus\mathcal F(k)$ and
$x_i(h)=0\ \forall i\in\mathcal R(k+1)$ and $h\geq k+1$.
All nodes in $\mathcal F(k)$ are detected from $\bm{x}(k)$ and
handled simultaneously, since a sequential order would change the outcome
when they share surviving neighbors.

Let $P(k)\in\mathbb{R}_{\geq0}^{n\times n}$ with entries $p_{ij}(k)$ be the
load-transfer matrix before failure detection at step $k$.  For an active
node $i$, $p_{ij}(k)$ is the fraction of its load sent to neighbor $j$,
and $p_{ii}(k)$ is the fraction retained at $i$.  Each row sums to one
because these fractions account for all of that node's load.  Thus, $P(k)$
is row stochastic. In summary, $G(k)$ specifies the available connections and
$P(k)$ specifies the transfer fractions along them.
At a failure-free step, the nominal load dynamics are
\begin{equation}
 \bm{x}(k+1)=P(k)^{\mathsf T}\bm{x}(k).
 \label{eq:nominal-phi}
\end{equation}
Let $\bm{1}\in\mathbb{R}^n$ denote the vector of all ones.  Since
$P(k)\bm{1}=\bm{1}$,
\begin{equation}
 \bm{1}^{\mathsf T}\bm{x}(k+1) = \bm{1}^{\mathsf T}P(k)^{\mathsf T}\bm{x}(k)
 =\bm{1}^{\mathsf T}\bm{x}(k),
 \label{eq:nominal-conservation}
\end{equation}
so nominal transfers preserve total load.
In order to construct the matrix $P$ and to assign transfer weights along the edges of a given network $G$, let $A\in\mathbb{R}_{\geq0}^{n\times n}$
with entries $a_{ij}$ be its weighted adjacency matrix,
$D=\operatorname{diag}(A\bm{1})$ the degree matrix with diagonal entries
$d_{ii}$, $L=D-A$, and $d_{\max}=\max_{i\in V} d_{ii}$.  For $d_{\max}>0$, the initial load-transfer matrix is
\begin{equation}
 P(0)=I-\varepsilon L,
 \qquad 0\leq\varepsilon\leq 1/d_{\max}.
 \label{eq:initial-transition}
\end{equation}
This diffusion-type update is used in consensus
dynamics~\cite{olfatisaber2007consensus}.  The bound on $\epsilon$ ensures nonnegative entries $p_{ij}(0)$, as well as the row-stochasticity, i.e., every row sums to one.  When $d_{\max}=0$, we have $P(0)=I$.

\subsection{Topology update at failure}

Let $\mathcal S(k)=\mathcal A(k)\setminus\mathcal F(k)$ be the surviving
nodes.  Their graph is $G(k+1)=G[\mathcal S(k)]$.
The matrix $P^{+}(k)$ retains dimension $n\times n$ associated with the initial node set $V$.
For $i\in\mathcal S(k)$, the entries read
\begin{subequations}
\let\quad\enspace 
\begin{numcases}{p^{+}_{ij}(k)=}
 0, & if $j\notin\mathcal S(k)$, \label{eq:p-zero-column}\\
 p_{ii}(k)+\sum_{\ell\in\mathcal F(k)}p_{i\ell}(k), & if $j\in\mathcal S(k)$, $j=i$, \label{eq:p-diagonal}\\
 p_{ij}(k), & if $j\in\mathcal S(k)$, $j\neq i$. \label{eq:p-survivors}
\end{numcases}
\end{subequations}
For every $i\in\mathcal F(k)\cup\mathcal R(k)$, the entries read
\begin{equation}
 p^{+}_{ij}(k)=
 \begin{cases}
  1, & \text{if } j=i,\\
  0, & \text{otherwise}.
 \end{cases}
 \label{eq:p-failed}
\end{equation}
By~\eqref{eq:p-zero-column} and~\eqref{eq:p-failed}, also
$p^{+}_{ji}(k)=0$ for $i\in\mathcal F(k)\cup\mathcal R(k)$,  $j\neq i$.
The intermediate load vector is computed as
\begin{equation}
 \widehat{\bm{x}}(k)=P^{+}(k)^{\mathsf T}\bm{x}(k).
 \label{eq:xhat}
\end{equation}
This vector gives the load distribution after the matrix update and before
the newly failed-node loads are reassigned.
By \eqref{eq:p-diagonal}, the transfer weight previously assigned
from a surviving node to a newly failed node is retained by that survivor as
self-weight, which is consistent with the weight handling of consensus iterations
over switching topologies~\cite{olfatisaber2007consensus}.  In contrast, naive row--column deletion without this reassignment can lead to unnecessary loss of load under the nominal dynamics.
For a three-node complete graph with $\varepsilon=1/3$,
failure of node~2 would give
\begin{equation}
 P^{\mathrm{naive}}=\frac{1}{3}
 \begin{bmatrix}1&0&1\\0&3&0\\1&0&1\end{bmatrix},\qquad
 P^+=\frac{1}{3}
 \begin{bmatrix}2&0&1\\0&3&0\\1&0&2\end{bmatrix}.
 \label{eq:three-node-update}
\end{equation}
The surviving rows of $P^{\mathrm{naive}}$ sum to $2/3$, so
nominal transfers loose $(x_1+x_3)/3$.  In $P^+$, diagonal compensation
preserves unit row sums and therefore also total load before redistribution.

\subsection{Failed-load redistribution and state update}
We now define how the load of each newly failed node is redistributed.
We assume that this load remains available for a one-time reassignment
using the node's connections immediately before removal.  Only surviving
neighbors can receive it and these incident connections are excluded from
subsequent nominal transfers.
For $i\in\mathcal F(k)$, define the active
receiving-neighbor set as
\begin{equation}
 \mathcal N_i(k)=\{j\in\mathcal S(k):(i,j)\in E\}.
 \label{eq:active-neighbors}
\end{equation}
Let
$q_{ij}(k)\geq0$ be the load transferred from failed node $i$ to
$j\in\mathcal N_i(k)$, and let $s_i(k)\geq0$ be its shed or unserved load.
All newly failed nodes satisfy the failed-load balance
\begin{equation}
 \sum_{j\in\mathcal N_i(k)}q_{ij}(k)+s_i(k)=x_i(k), \quad \text{for all}
 i\in\mathcal F(k).
 \label{eq:failed-balance}
\end{equation}
Thus, every part of newly failed-node load is either assigned to an admissible
receiver or marked explicitly as unserved.

At each time step, the state is updated in the following order.  First,
$\mathcal F(k)$ is derived from $\bm{x}(k)$ with \eqref{eq:failure-set}. The load-transfer matrix is then updated to $P^{+}(k)$, and the
intermediate load vector $\hat{x}(k)$ is computed with~\eqref{eq:xhat}. Then a redistribution scheme
determines $q_{ij}(k)$ and $s_i(k)$.  For each $j\in\mathcal S(k)$, let
$\mathcal L_j(k)=\{i\in\mathcal F(k)\mid j\in\mathcal N_i(k)\}$
denote the set of newly failed nodes that can transfer load to $j.$ Then the state update is
\begin{subequations}
\begin{align}
 x_j(k+1)&=\widehat{x}_j(k)
 +\sum_{i\in\mathcal L_j(k)}q_{ij}(k), \quad
 \text{for } j\in\mathcal S(k),                               \label{eq:survivor-state}\\
 x_i(k+1)&=0,\quad \text{for }
 i\in\mathcal F(k)\cup\mathcal R(k).                   \label{eq:failed-state}
\end{align}
\end{subequations}
Equation~\eqref{eq:survivor-state} dispatches the load at each surviving node,
whereas~\eqref{eq:failed-state} sets all failed-node states to zero. Then the update $\mathcal{A}(k+1) = \mathcal{A}(k)\setminus \mathcal{F}(k)$ and $P(k+1) = P^+(k)$ is performed. This update order is used
throughout the analyses and simulations.

\begin{remark}
The update~\eqref{eq:survivor-state}--\eqref{eq:failed-state} can be equivalently written as
$\bm{x}(k+1)=P^{+}(k)^{\mathsf T}\bm{x}(k)+\bm{\phi}(k)$, where
$\bm{\phi}(k)\in\mathbb R^n$ accounts for failed-load removals and
redistribution. Its entries are
\begin{equation}
 \phi_j(k)=
 \begin{cases}
  \displaystyle\sum_{i\in\mathcal L_j(k)}q_{ij}(k),
   & \text{if } j\in\mathcal S(k),\\[2ex]
  -x_j(k), &  \text{if } j\in\mathcal F(k),\\[0.5ex]
  0, & \text{if } j\in\mathcal R(k).
 \end{cases}
 \label{eq:phi-meaning}
\end{equation}
At a failure-free step, $\mathcal F(k)=\varnothing$ implies
$P^{+}(k)=P(k)$ and $\bm{\phi}(k)=\bm{0}$, so the full update
reduces to~\eqref{eq:nominal-phi}.
\end{remark}

\section{Load Balance and Lossless Redistribution}
\label{sec:analysis}
In this section, we first show that the proposed matrix update preserves row stochasticity and derive the resulting total-load balance.  We then characterize when the loads
of simultaneously failed nodes can be fully redistributed without exceeding
the capacities of the surviving nodes.

We start by showing that \eqref{eq:p-diagonal}
preserves the stochastic structure required by
\eqref{eq:nominal-conservation}.

\begin{proposition}
\label{prop:stochasticity}
If $P(k)$ is nonnegative and row stochastic, then the update
\eqref{eq:p-zero-column}--\eqref{eq:p-failed} produces a nonnegative,
row-stochastic $P^{+}(k)$.  If $P(k)$ is also symmetric, then $P^{+}(k)$ is
symmetric and doubly stochastic.
\end{proposition}

\begin{proof}
Nonnegativity follows because the update only retains nonnegative entries,
sets entries to zero or one, and adds removed nonnegative weights to 
diagonal.  For entries $i\in\mathcal S(k)$, summing the updated rows gives
\begin{equation}
\csname ctagsplit@false\endcsname
\begin{split}
\sum_{j\in V}p^+_{ij}(k)
 &=p^+_{ii}(k)
   +\sum_{j\in\mathcal S(k)\setminus\{i\}}p^+_{ij}(k)
   +\underbrace{\sum_{j\notin\mathcal S(k)}
     p^+_{ij}(k)}_{=0\text{ by }\eqref{eq:p-zero-column}}\\
 &=\underbrace{p_{ii}(k)+\sum_{\ell\in\mathcal F(k)}p_{i\ell}(k)}_{\text{by }\eqref{eq:p-diagonal}}
   +\underbrace{\sum_{j\in\mathcal S(k)\setminus\{i\}}p_{ij}(k)}_{\text{by }\eqref{eq:p-survivors}}\\
 &=\sum_{j\in \mathcal{F}(k) \cup \mathcal{S}(k)}p_{ij}(k) = \sum_{j\in V}p_{ij}(k)=1,\raisetag{10pt}
\end{split}
\label{eq:prop1-rowsum}
\end{equation}
since entries $p_{ij}(k)$ associated with previously failed nodes $j \in \mathcal{R}(k)$ are already zero outside
their diagonals.  For $i\in\mathcal F(k)\cup\mathcal R(k)$, the updated row is
an identity row and therefore also sums to one.  If $P(k)$ is symmetric, all
surviving off-diagonal pairs remain unchanged and both entries incident to a
failed node are set to zero.  Hence $P^{+}(k)$ remains symmetric. Symmetry and
unit row sums then imply unit column sums.
\end{proof}

Thus, the updated matrix preserves total load under the nominal dynamics,
whether one or several nodes fail at the same time. The next proposition
gives the resulting total-load balance.
\begin{proposition}
\label{prop:load-balance}
In the state update~\eqref{eq:survivor-state}--\eqref{eq:failed-state},
suppose that redistribution flows and unserved loads satisfy
\eqref{eq:failed-balance}.  Then,
\begin{equation}
 \bm{1}^{\mathsf T}\bm{x}(k+1)
 =\bm{1}^{\mathsf T}\bm{x}(k)
 -\sum_{i\in\mathcal F(k)}s_i(k).
 \label{eq:total-load-balance}
\end{equation}
\end{proposition}

\begin{proof}
By Proposition~\ref{prop:stochasticity},
$P^+(k)\bm{1}=\bm{1}.$ So, by~\eqref{eq:xhat},
{\relpenalty=10000\binoppenalty=10000
$\bm{1}^{\mathsf T}\widehat{\bm{x}}(k)
=\bm{1}^{\mathsf T}P^+(k)^{\mathsf T}\bm{x}(k)
=\bigl(P^+(k)\bm{1}\bigr)^{\mathsf T}\bm{x}(k)
=\bm{1}^{\mathsf T}\bm{x}(k)$}.
The isolated column of each newly failed node
$i\in\mathcal F(k)$ gives
\begin{equation}
 \widehat{x}_i(k)
 =\sum_{j\in V}p^+_{ji}(k)x_j(k)
 \overset{\eqref{eq:p-zero-column},\eqref{eq:p-failed}}{=}x_i(k).
 \label{eq:xhat-failed}
\end{equation}
Previously failed nodes $i\in\mathcal R(k)$ have zero load by the irreversible-failure
rule, and their isolated identity form also gives $\widehat{x}_i(k)=0.$  Since $\mathcal S(k)$, $\mathcal F(k)$, and
$\mathcal R(k)$ partition $V$, the intermediate load on surviving nodes is
\begin{equation}
\begin{split}
 \sum_{j\in\mathcal S(k)}\widehat{x}_j(k)
 &=\bm{1}^{\mathsf T}\widehat{\bm{x}}(k)
   -\sum_{i\in\mathcal F(k)}\widehat{x}_i(k) - \sum_{i\in\mathcal R(k)}\widehat{x}_i(k)\\
 &=\bm{1}^{\mathsf T}\bm{x}(k)
   -\sum_{i\in\mathcal F(k)}x_i(k).
\end{split}
\label{eq:surviving-intermediate}
\end{equation}
By \eqref{eq:failed-state}, we have that
\begin{equation}
\begin{split}
 \bm{1}^{\mathsf T}\bm{x}(k+1)
 & ~=\sum_{j\in\mathcal S(k)}x_j(k+1)\\
 &\overset{\eqref{eq:survivor-state}}{=}\sum_{j\in\mathcal S(k)}\widehat{x}_j(k)
   +\sum_{j\in\mathcal S(k)}\sum_{i\in\mathcal L_j(k)}q_{ij}(k)\\
 &\overset{\phantom{\eqref{eq:survivor-state}}}{=}\sum_{j\in\mathcal S(k)}\widehat{x}_j(k)
   +\sum_{i\in\mathcal F(k)}\sum_{j\in\mathcal N_i(k)}q_{ij}(k)\\
 &\overset{\eqref{eq:surviving-intermediate},\eqref{eq:failed-balance}}{=}\bm{1}^{\mathsf T}\bm{x}(k)
   -\sum_{i\in\mathcal F(k)}x_i(k)\\
 &\qquad\qquad \qquad
   +\sum_{i\in\mathcal F(k)}\bigl(x_i(k)-s_i(k)\bigr)\\
 &=\bm{1}^{\mathsf T}\bm{x}(k)
   -\sum_{i\in\mathcal F(k)}s_i(k).
\end{split}
\label{eq:prop2-chain}
\end{equation}
This completes the proof.
\end{proof}

If $s_i(k)=0$ for every $i\in\mathcal F(k)$, then
\eqref{eq:total-load-balance} gives lossless redistribution. Otherwise, total load decreases by $\sum_{i\in\mathcal F(k)}s_i(k)$, which is precisely the amount that could
not be redistributed.

We next characterize whether all newly failed loads can be transferred without
immediately overloading a survivor.  Define the residual capacity after the
transfers $\widehat{\bm x}(k)=P^{+}(k)^{\mathsf T}\bm x(k)$ by
\begin{equation}
 r_j(k)=[c_j-\widehat{x}_j(k)]_+.
 \label{eq:residual-capacity}
\end{equation}
where $[z]_+=\max\{z,0\}$. For $\mathcal U\subseteq\mathcal F(k)$, let
\begin{equation}
 \mathcal N(\mathcal U,k)=\bigcup_{i\in\mathcal U}\mathcal N_i(k)
 \label{eq:receiver-union}
\end{equation}
denote the set of active receivers adjacent to at least one node in
$\mathcal U$.
\begin{theorem}
\label{thm:simultaneous-containment}
Let $\mathcal F(k)$ be the set of nodes that fail simultaneously at
time~$k$, and assume that $\widehat{x}_j(k)\leq c_j$ for every $j\in\mathcal S(k)$.  A
redistribution with nonnegative flows only to $j\in\mathcal N_i(k)$,
satisfying~\eqref{eq:failed-balance} with $s_i(k)=0$ and
$x_j(k+1)\leq c_j$ for every $j\in\mathcal S(k)$, exists if and only if
\begin{equation}
 \sum_{i\in\mathcal U}x_i(k)
 \leq
 \sum_{j\in\mathcal N(\mathcal U,k)}r_j(k)
 \quad\text{for every }\mathcal U\subseteq\mathcal F(k).
 \label{eq:simultaneous-containment}
\end{equation}
Moreover, the minimum unavoidable unserved load is
\begin{equation}
 \sigma^\star(k)=
 \max_{\mathcal U\subseteq\mathcal F(k)}
 \left[
  \sum_{i\in\mathcal U}x_i(k)
  -\sum_{j\in\mathcal N(\mathcal U,k)}r_j(k)
 \right]_+.
 \label{eq:cut-deficiency}
\end{equation}
\end{theorem}

\begin{figure}[!t]
 \centering
 \scalebox{0.85}{\begin{tikzpicture}[
  every node/.style={font=\small},
  vertex/.style={circle,draw,fill=black!5,minimum size=6mm,inner sep=1pt},
  terminal/.style={circle,fill,minimum size=2mm,inner sep=0pt},
  capacity/.style={font=\small,fill=white,inner sep=1pt}
]
 \node[terminal,label=left:source] (source) at (0,0) {};
 \node[vertex] (f1) at (1.6,0.7) {$i_1$};
 \node[vertex] (f2) at (1.6,-0.7) {$i_2$};
 \node[vertex] (r1) at (4,0.7) {$j_1$};
 \node[vertex] (r2) at (4,-0.7) {$j_2$};
 \node[terminal,label=right:sink] (sink) at (5.6,0) {};
 \node at (1.6,1.25) {Failed nodes};
 \node at (4,1.25) {Receivers};
 \draw[-latex] (source) -- node[capacity,above] {$x_{i_1}(k)$} (f1);
 \draw[-latex] (source) -- node[capacity,below] {$x_{i_2}(k)$} (f2);
 \draw[-latex] (f1) -- node[capacity,above] {$M$} (r1);
 \draw[-latex] (f2) -- node[capacity,above] {$M$} (r1);
 \draw[-latex] (f2) -- node[capacity,below] {$M$} (r2);
 \draw[-latex] (r1) -- node[capacity,above] {$r_{j_1}(k)$} (sink);
 \draw[-latex] (r2) -- node[capacity,below] {$r_{j_2}(k)$} (sink);
\end{tikzpicture}}
 \caption{Auxiliary flow network for the connections in
 Example~\ref{ex:shared-receiver}. Edge labels indicate the corresponding capacities.}
 \label{fig:redistribution-flow-network}
\end{figure}
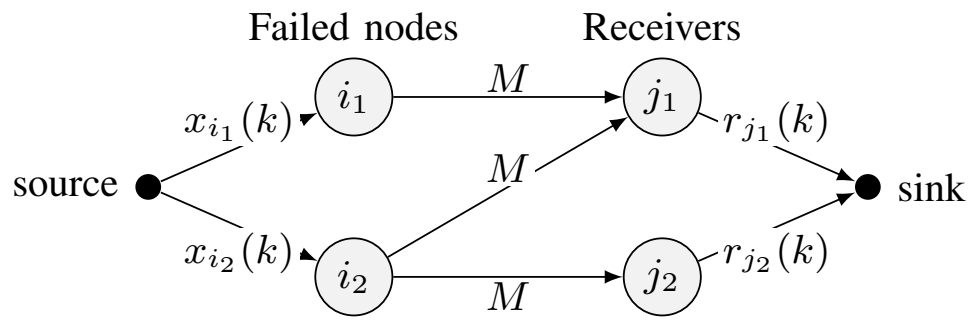

\begin{proof}
If $\mathcal F(k)=\varnothing$, both claims hold with $\sigma^\star(k)=0$.
Otherwise, let $M=\sum_{i\in\mathcal F(k)}x_i(k)>0$.  The idea is to
represent the redistribution as a flow from the failed nodes to the
receivers and to apply the max-flow/min-cut theorem~\cite{ford1956maximal}.

First, we construct an auxiliary directed network as illustrated in
Fig.~\ref{fig:redistribution-flow-network}.  To distinguish it from $G$, we
refer to its elements as vertices and arcs.  Connect the source to each
$i\in\mathcal F(k)$ with capacity $x_i(k)$, vertex $i$ to each
$j\in\mathcal N_i(k)$ with capacity $M$, and each receiver
$j\in\mathcal S(k)$ to the sink with capacity $r_j(k)$.  Since
$M\geq x_i(k)$, the failed-to-receiver arcs never restrict a transfer.

Then, we relate flows in this network to redistributions.  A feasible flow
is nonnegative, respects arc capacities, and conserves flow at every vertex
except the source and sink.  Its value is the total flow leaving the source.
Let the flow on arc $(i,j)$ be $q_{ij}(k)$.  Conservation at $i$ makes the
source-to-$i$ flow equal to $\sum_{j\in\mathcal N_i(k)}q_{ij}(k)$.
Defining $s_i(k)$ as the unused capacity of this arc, i.e., $x_i(k)$ minus
this flow, gives the failed-load balance~\eqref{eq:failed-balance}.  Conservation at $j$
makes the $j$-to-sink flow equal to $\sum_{i\in\mathcal L_j(k)}q_{ij}(k)$.
Since $\widehat{x}_j(k)\leq c_j$, we have $r_j(k)=c_j-\widehat{x}_j(k)$, so
by~\eqref{eq:survivor-state} the capacity of this arc is respected if and
only if $x_j(k+1)\leq c_j$.  Thus, feasible flows correspond to admissible
redistributions within the receiver capacities.  A flow of value $M$
saturates every source arc and is equivalent to lossless redistribution.

Finally, we determine the maximum flow value.  A source--sink cut partitions
the vertices with the source and sink on opposite sides. Its capacity is
the sum of capacities of the arcs, directed from the
source side to the sink side crossing it.  The max-flow/min-cut theorem equates the maximum flow value and the
minimum cut capacity.  The cut with only the source on its side has
capacity $M$, and any cut crossing a failed-to-receiver arc has capacity
at least $M$.  Consider a cut crossing no such arc. Let $\mathcal U$ contain the
failed vertices on the source side.  Every receiver in
$\mathcal N(\mathcal U,k)$ must then also lie on that side. For fixed
$\mathcal U$, placing only these receivers on the source side minimizes
the cut capacity. The crossing arcs run from the source to failed vertices
in $\mathcal F(k)\setminus\mathcal U$ and from receivers in
$\mathcal N(\mathcal U,k)$ to the sink, giving
\begin{align}
 C(\mathcal U)&=\sum_{i\in\mathcal F(k)\setminus\mathcal U}x_i(k)
 +\sum_{j\in\mathcal N(\mathcal U,k)}r_j(k)\nonumber\\
 &=M-\sum_{i\in\mathcal U}x_i(k)
 +\sum_{j\in\mathcal N(\mathcal U,k)}r_j(k).
 \label{eq:cut-capacity}
\end{align}
Hence, the maximum flow value is
$\min_{\mathcal U\subseteq\mathcal F(k)}C(\mathcal U)$, which equals $M$
exactly when~\eqref{eq:simultaneous-containment} holds. By~\eqref{eq:cut-capacity},
the minimum unserved load is therefore
\[
\begin{aligned}
 \sigma^\star(k)
 &=M-\min_{\mathcal U\subseteq\mathcal F(k)}C(\mathcal U)\\
 &=\max_{\mathcal U\subseteq\mathcal F(k)}
 \left[
  \sum_{i\in\mathcal U}x_i(k)
  -\sum_{j\in\mathcal N(\mathcal U,k)}r_j(k)
 \right]_+.
\end{aligned}
\]
\end{proof}

The following example shows why comparing total failed load with total
residual capacity is insufficient when the admissible transfers are restricted by the network topology.
\begin{example}
\label{ex:shared-receiver}
Consider two simultaneously failed nodes $i_1$ and $i_2$, each carrying
$0.8$ units, and two active receivers $j_1$ and $j_2$ with residual
capacities $r_{j_1}=\rho$ and $r_{j_2}=1$.  Let $i_1$ be connected only to
$j_1$, while $i_2$ is connected to both receivers.  The aggregate residual
capacity is sufficient when $\rho\geq0.6$.  However, applying
\eqref{eq:simultaneous-containment} to $\mathcal U=\{i_1\}$ gives the tighter
condition $\rho\geq0.8$.  Moreover, equal redistribution sends $0.8$ from
$i_1$ to $j_1$ and $0.4$ from $i_2$ to $j_1$, so it violates that receiver's capacity
unless $\rho\geq1.2$.  Hence, aggregate residual capacity does not capture
the topology-induced bottleneck.
\end{example}

\begin{remark}[Single-node case]
Under the assumptions of
Theorem~\ref{thm:simultaneous-containment}, for $\mathcal F(k)=\{i\}$,
lossless redistribution within receiver capacities exists if and only if
$x_i(k)\leq\sum_{j\in\mathcal N_i(k)}r_j(k)$.
\end{remark}

\section{Optimization-Based Redistribution Control}
\label{sec:control}

This section presents a two-stage linear programming approach for
redistributing newly failed-node loads.  The first stage minimizes unserved
load, while the second stage minimizes the maximum capacity
utilization among surviving nodes. As opposed to equal redistribution, we consider residual capacities and coordinate
failed nodes sharing a receiver by actively selecting the redistribution flows.
Let $\bm q(k)$ collect $q_{ij}(k)$ for $i\in\mathcal F(k)$ and
$j\in\mathcal N_i(k)$, and let $\bm s(k)$ collect $s_i(k)$ for
$i\in\mathcal F(k)$.

\subsection{Stage 1: Minimum shedding}

The first stage finds the minimum unserved load $\sigma^\star(k)$ consistent with
the surviving-node capacities.

\begin{problem}[Minimum-shedding redistribution]
\label{prob:stage-one}
When $\mathcal F(k)\neq\varnothing$, let us compute $\sigma^\star(k)$ as follows:
\begin{subequations}
\label{eq:stage-one}
\begin{align}
 \underset{\bm q(k),\,\bm s(k)}{\operatorname{minimize}}\quad
 &\sum_{i\in\mathcal F(k)}s_i(k)
 \label{eq:stage-one-objective}
 \shortintertext{subject to}
 \widehat{\bm{x}}(k)&=P^{+}(k)^{\mathsf T}\bm{x}(k),
 \label{eq:stage-one-nominal}\\
 \widehat{x}_j(k)+\sum_{i\in\mathcal L_j(k)}q_{ij}(k)&\leq c_j,\quad
 j\in\mathcal S(k),                                    \label{eq:stage-one-capacity}\\
 \sum_{j\in\mathcal N_i(k)}q_{ij}(k)+s_i(k)&=x_i(k),\quad
 i\in\mathcal F(k),                                    \label{eq:stage-one-balance}\\
 \bm q(k)&\geq\bm 0,\quad \bm s(k)\geq\bm 0.
 \label{eq:stage-one-positive}
\end{align}
\end{subequations}
\end{problem}
Constraint~\eqref{eq:stage-one-capacity} enforces capacity constraints, whereas
\eqref{eq:stage-one-balance} permits load to leave only through $s_i(k)$.
Under the nominal survivor-feasibility assumption
$\widehat{x}_j(k)\leq c_j$ of
Theorem~\ref{thm:simultaneous-containment}, Problem~\ref{prob:stage-one} is
always feasible for $q_{ij}(k)=0$ and $s_i(k)=x_i(k)$.  If
$\widehat{x}_j(k)>c_j$ for any surviving node before redistribution,
Problem~\ref{prob:stage-one} is infeasible because the nonnegative
redistribution flows $q_{ij}(k)$ cannot reduce that node's existing load.  No capacity
guarantee is made in this case. 
Note, that problem~1 leads to $\sigma^\star(k)=0$ if lossless redistribution is feasible, and
$\sigma^\star(k)>0$ if a minimum load that must remain unserved exists.

\subsection{Stage 2: Utilization balancing}
Stage~2 minimizes maximum capacity utilization while allowing a total unserved
load of at most $\sigma^\star(k)+\delta_{\mathrm{LP}}$, where
$\delta_{\mathrm{LP}}>0$ is a small tolerance that accounts for the numerical
precision of the solver.  This maximizes the smallest
unused capacity fraction $1-x_j(k+1)/c_j$.

\begin{problem}[Utilization-balanced redistribution]
\label{prob:stage-two}
Given $\sigma^\star(k)$ from Problem~\ref{prob:stage-one}, solve
\begin{subequations}
\label{eq:stage-two}
\begin{align}
 \underset{\bm q(k),\,\bm s(k),\,\eta}{\operatorname{minimize}}\quad&\eta
 \label{eq:stage-two-objective}
 \shortintertext{subject to}
 \widehat{\bm{x}}(k)&=P^{+}(k)^{\mathsf T}\bm{x}(k),
 \label{eq:stage-two-nominal}\\
 \widehat{x}_j(k)+\sum_{i\in\mathcal L_j(k)}q_{ij}(k)&\leq \eta c_j,\quad
 j\in\mathcal S(k),                                    \label{eq:stage-two-utilization}\\
 \sum_{j\in\mathcal N_i(k)}q_{ij}(k)+s_i(k)&=x_i(k),\quad
 i\in\mathcal F(k),                                    \label{eq:stage-two-balance}\\
 \sum_{i\in\mathcal F(k)}s_i(k)&\leq \sigma^\star(k)+\delta_{\mathrm{LP}},
 \label{eq:stage-two-shed}\\
 \bm q(k)&\geq\bm 0,\quad \bm s(k)\geq\bm 0,\quad\eta\geq0.
 \label{eq:stage-two-positive}
\end{align}
\end{subequations}
\end{problem}
Since $c_j>0$, $\eta$ bounds the largest normalized surviving-node load
after redistribution.  Problem~\ref{prob:stage-two} is feasible because any
optimizer of Problem~\ref{prob:stage-one}, together with $\eta=1$, satisfies
all its constraints.  Let $\eta^\star(k)$ denote the optimal value of
Problem~\ref{prob:stage-two}.  Hence, $\eta^\star(k)\leq1$, and every surviving node
remains within its capacity.  The Stage~2 solution is therefore also feasible for
Stage~1. Stage~2 provides
\begin{equation}
 \sigma^\star(k)\leq\sum_{i\in\mathcal F(k)}s_i(k)
 \leq\sigma^\star(k)+\delta_{\mathrm{LP}}.
 \label{eq:sigma-opt}
\end{equation}
Thus, balancing receiver utilization increases the unserved load by at most
the desired numerical tolerance. When $\mathcal F(k)\neq\varnothing$, Problem~\ref{prob:stage-one} is solved
first, and Problem~\ref{prob:stage-two} is then solved using
$\sigma^\star(k)$.  Both problems concern only the current redistribution.

For subsequent nominal updates, let
$B=[p^+_{ij}(k)]_{i,j\in\mathcal S(k)}$ be the surviving-node block, and let
$\bm c_{\mathcal S}$ and $\bm x_{\mathcal S}(h)$ denote the capacity and load
vectors only containing the values of $\mathcal S(k)$.

\begin{proposition}[Failure-free invariance]
\label{cor:failure-free-invariance}
Suppose $\mathcal S(k)\neq\varnothing$, $B\geq0$ is row stochastic, and
$\bm c_{\mathcal S}>0$.  With $B$ fixed and no further external input, the
nominal updates $\bm x_{\mathcal S}(h+1)=B^{\mathsf T}\bm x_{\mathcal S}(h)$
preserve $0\leq\bm x_{\mathcal S}(h)\leq\bm c_{\mathcal S}$ for every initial
state within these bounds if and only if
\begin{equation}
 B^{\mathsf T}\bm c_{\mathcal S}=\bm c_{\mathcal S}.
\end{equation}
If the condition holds and Problem~\ref{prob:stage-one} is feasible, applying
the Stage~2 solution gives no additional failures at any $h\geq k+1$.
\end{proposition}

\begin{proof}
Nonnegativity of $B$ gives
\[
 0\leq B^{\mathsf T}\bm x_{\mathcal S}(h)
 \leq B^{\mathsf T}\bm c_{\mathcal S}=\bm c_{\mathcal S},
\]
so the condition is sufficient by induction.  Necessity follows by taking
the feasible state $\bm x_{\mathcal S}(h)=\bm c_{\mathcal S}$, which requires
$B^{\mathsf T}\bm c_{\mathcal S}\leq\bm c_{\mathcal S}$.  Since
$B\bm 1=\bm 1$, both vectors have the same sum. Since no entry on the left can be larger, none can be smaller either, so equality holds
componentwise. Finally, the bound $\eta^\star(k)\leq1$ resulting from
Problem~\ref{prob:stage-two} places the state after redistribution within
these invariant bounds.
\end{proof}

Uniform capacities with a doubly stochastic $B$, as in Section~\ref{sec:results},
satisfy this condition, whereas for unequal capacities a redistribution within
all capacities at step~$k$ can still lead to later overloads.  Under this
condition, the surviving nodes also satisfy
$\widehat{\bm x}_{\mathcal S}(k)=B^{\mathsf T}\bm x_{\mathcal S}(k)
\leq\bm c_{\mathcal S}$, so Problem~\ref{prob:stage-one} is always feasible.

\begin{remark}[Directed networks]
\label{rem:directed}
Propositions~\ref{prop:stochasticity} and~\ref{prop:load-balance},
Theorem~\ref{thm:simultaneous-containment}, and
Problems~\ref{prob:stage-one} and~\ref{prob:stage-two} do not rely on
symmetry. They also hold for a directed graph $G$ in which an edge $(i,j)$
only allows load transfer from $i$ to $j$, if $\mathcal N_i(k)$ denotes the
set of surviving out-neighbors of node $i$ and $P(0)=I-\varepsilon L$ is
constructed with the out-degree Laplacian $L=D_{\mathrm{out}}-A$ and
$0\leq\varepsilon\leq1/d^{\mathrm{out}}_{\max}$. In this case, however, $B$
is in general not doubly stochastic. Hence, the condition of
Proposition~\ref{cor:failure-free-invariance} may fail even for uniform
capacities, and a redistribution within all capacities at step~$k$ can still
lead to later overloads.
\end{remark}

\section{Numerical Evaluation}
\label{sec:results}

This section compares equal redistribution, residual-capacity redistribution, and controlled redistribution on grid and scale-free networks.  The evaluation focuses on secondary
failures and unserved load.

\subsection{Simulation setup and metrics}

We consider the following networks: $10\times10$ grid and scale-free networks with $n=100$, $c_i=1\ \forall i \in V$,
baseline load $x_i=0.1\ \forall i \in V$, and
$P(0)$ as in~\eqref{eq:initial-transition} with fixed $\varepsilon=0.5/d_{\max}$.
In each scenario, we consider two nodes, with a common surviving neighbor. Both of them receive an
additional load $\Delta\in\{1,1.25,2,3,4.5\}$ causing them to fail simultaneously.  Their
common surviving neighboring node ensures that the two failed nodes
compete for at least one receiver.  On the $10\times10$ grid, we consider all
$322$ such node pairs, which gives $1610$ scenarios per redistribution method over all values of $\Delta.$
For each of $30$ Barabasi--Albert graphs~\cite{barabasi1999emergence} with
attachment parameter $2$, we sample $20$ such pairs uniformly without
replacement, which gives $3000$ scenarios per method.

We compare equal redistribution, residual-capacity proportional
redistribution, and the controlled redistribution of
Section~\ref{sec:control}.  The residual-capacity method distributes each
failed-node load among its active receivers in proportion to the residual
capacities $r_j(k)$ defined in~\eqref{eq:residual-capacity}.  If all
corresponding residual capacities are zero, it uses equal redistribution.
The optimization problems associated with controlled redistribution are solved with
\texttt{scipy.optimize.linprog} using HiGHS.  Let $\mathcal F_{0}$ be
the set of initially disturbed nodes, and let $\mathcal F_{\mathrm e}$
contain all nodes that fail during the cascade.  The secondary-failure
fraction is
$f_{\mathrm{sec}}=|\mathcal F_{\mathrm e}\setminus\mathcal F_{0}|/
(n-|\mathcal F_{0}|)$.
The unserved-load fraction is defined as the total unserved load accumulated over the
cascade divided by the initial total load.
For the uniform-capacity simulations,
Proposition~\ref{cor:failure-free-invariance} ensures that, once a step
produces no new failure, subsequent nominal updates cannot cause another one.
All runs pass the state, flow,
stochasticity, and balance checks, and no solver fails or iteration limit is
reached.

\subsection{Simulation results and discussion}

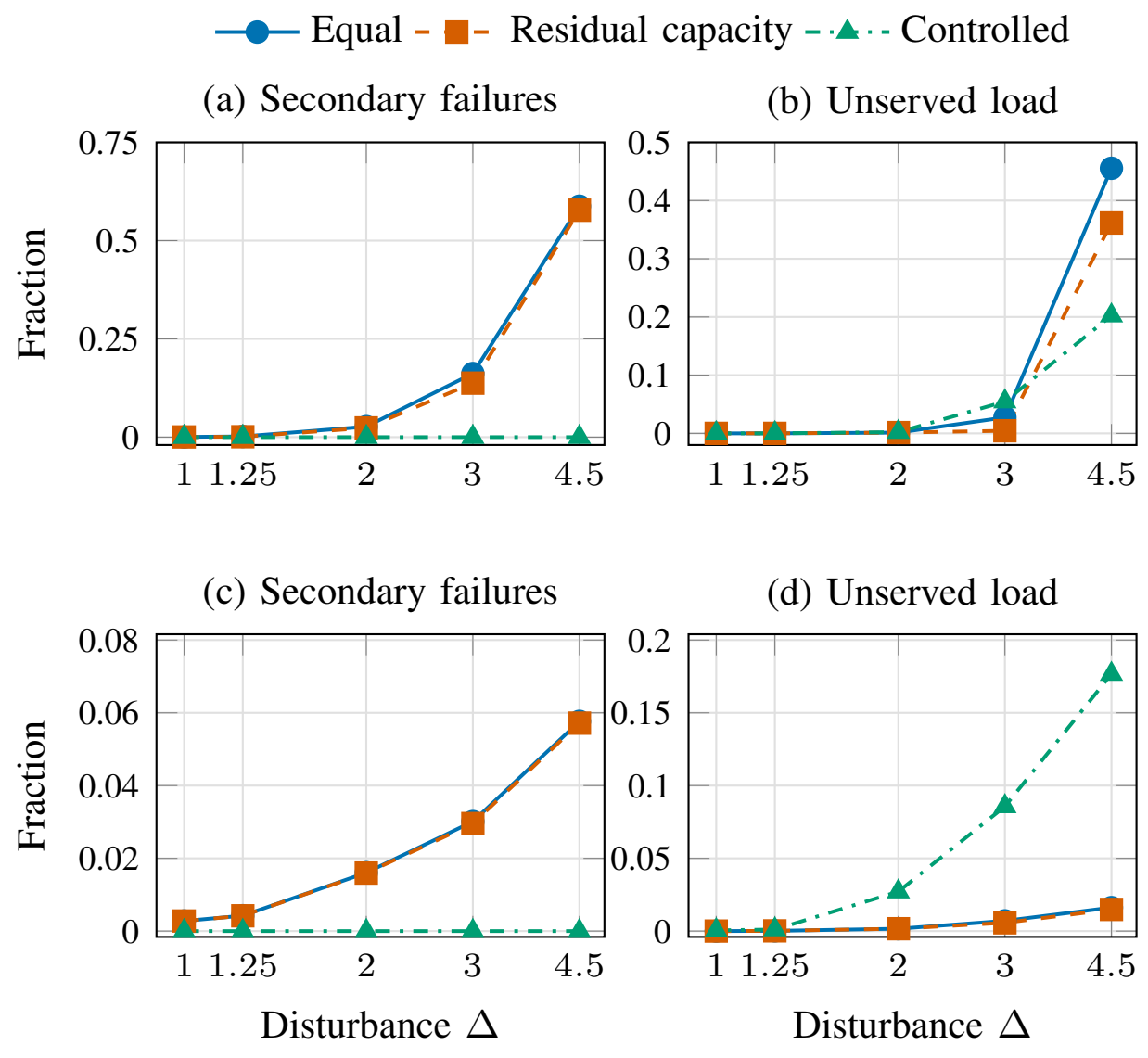
\begin{figure}[!t]
 \centering
 \begin{tikzpicture}
\definecolor{cfEqual}{HTML}{0072B2}
\definecolor{cfResidual}{HTML}{D55E00}
\definecolor{cfControlled}{HTML}{009E73}
\begin{groupplot}[
  group style={group size=2 by 2, horizontal sep=0.60cm, vertical sep=1.35cm},
  scale only axis, width=0.37\linewidth, height=0.25\linewidth,
  xmode=log, log ticks with fixed point,
  xmin=0.9, xmax=4.95, ymin=-0.02, ymax=1.02,
  xtick={1,1.25,2,3,4.5},
  ytick={0,0.25,0.5,0.75,1},
  xticklabel style={font=\scriptsize},
  tick label style={font=\scriptsize}, label style={font=\footnotesize},
  every axis title shift=0pt,
  title style={font=\footnotesize, at={(0.5,1)}, anchor=south, yshift=1pt},
  grid=major, grid style={black!12},
  every axis plot/.append style={line width=0.75pt},
  legend style={font=\footnotesize, draw=none},
]
\nextgroupplot[title={(a) Secondary failures}, ylabel={Fraction}, ymin=-0.02, ymax=0.75, ytick={0,0.25,0.5,0.75}, scaled y ticks=false, yticklabels={0,0.25,0.5,0.75}, legend to name={simultaneouscombinedlegend}, legend columns=-1]
  \addlegendimage{cfEqual, solid, mark=*, mark options={solid, draw=cfEqual, fill=cfEqual}}
  \addlegendentry{Equal}
  \addlegendimage{cfResidual, dashed, mark=square*, mark options={solid, draw=cfResidual, fill=cfResidual}}
  \addlegendentry{Residual capacity}
  \addlegendimage{cfControlled, dashdotted, mark=triangle*, mark options={solid, draw=cfControlled, fill=cfControlled}}
  \addlegendentry{Controlled}
  \addplot+[cfEqual, solid, mark=*, mark options={solid, draw=cfEqual, fill=cfEqual}] coordinates {(1,0.00050703511) (1.25,0.0017746229) (2,0.026746102) (3,0.16187096) (4.5,0.5876537)};
  \addplot+[cfResidual, dashed, mark=square*, mark options={solid, draw=cfResidual, fill=cfResidual}] coordinates {(1,0.00050703511) (1.25,0.0010140702) (2,0.023450374) (3,0.13740652) (4.5,0.57738623)};
  \addplot+[cfControlled, dashdotted, mark=triangle*, mark options={solid, draw=cfControlled, fill=cfControlled}] coordinates {(1,0) (1.25,0) (2,0) (3,0) (4.5,0)};
\nextgroupplot[title={(b) Unserved load}, ymin=-0.02, ymax=0.5, ytick={0,0.1,0.2,0.3,0.4,0.5}, scaled y ticks=false, yticklabels={0,0.1,0.2,0.3,0.4,0.5}]
  \addplot+[cfEqual, solid, mark=*, mark options={solid, draw=cfEqual, fill=cfEqual}] coordinates {(1,0) (1.25,0) (2,0.0013309672) (3,0.027401017) (4.5,0.45542698)};
  \addplot+[cfResidual, dashed, mark=square*, mark options={solid, draw=cfResidual, fill=cfResidual}] coordinates {(1,0) (1.25,0) (2,0.0013309672) (3,0.0044720837) (4.5,0.36148974)};
  \addplot+[cfControlled, dashdotted, mark=triangle*, mark options={solid, draw=cfControlled, fill=cfControlled}] coordinates {(1,8.3333333e-11) (1.25,8e-11) (2,0.0021295475) (3,0.054658385) (4.5,0.20235371)};
\nextgroupplot[title={(c) Secondary failures}, xlabel={Disturbance $\Delta$}, ylabel={Fraction}, ymin=-0.0016, ymax=0.0816, ytick={0,0.02,0.04,0.06,0.08}, scaled y ticks=false, yticklabels={0,0.02,0.04,0.06,0.08}]
  \addplot+[cfEqual, solid, mark=*, mark options={solid, draw=cfEqual, fill=cfEqual}] coordinates {(1,0.0027721088) (1.25,0.0042346939) (2,0.016054422) (3,0.030102041) (4.5,0.057636054)};
  \addplot+[cfResidual, dashed, mark=square*, mark options={solid, draw=cfResidual, fill=cfResidual}] coordinates {(1,0.0027721088) (1.25,0.0042006803) (2,0.015952381) (3,0.029557823) (4.5,0.057159864)};
  \addplot+[cfControlled, dashdotted, mark=triangle*, mark options={solid, draw=cfControlled, fill=cfControlled}] coordinates {(1,0) (1.25,0) (2,0) (3,0) (4.5,0)};
\nextgroupplot[title={(d) Unserved load}, xlabel={Disturbance $\Delta$}, ymin=-0.004, ymax=0.204, ytick={0,0.05,0.1,0.15,0.2}, scaled y ticks=false, yticklabels={0,0.05,0.1,0.15,0.2}]
  \addplot+[cfEqual, solid, mark=*, mark options={solid, draw=cfEqual, fill=cfEqual}] coordinates {(1,0) (1.25,0.00014833333) (2,0.0016171627) (3,0.007038457) (4.5,0.01615745)};
  \addplot+[cfResidual, dashed, mark=square*, mark options={solid, draw=cfResidual, fill=cfResidual}] coordinates {(1,0) (1.25,0.00014833333) (2,0.0016179328) (3,0.0056914876) (4.5,0.014877316)};
  \addplot+[cfControlled, dashdotted, mark=triangle*, mark options={solid, draw=cfControlled, fill=cfControlled}] coordinates {(1,0.00050000008) (1.25,0.0010800001) (2,0.027142857) (3,0.0856875) (4.5,0.17664912)};
\end{groupplot}
\node[above=4mm] at ($(group c1r1.north west)!0.5!(group c2r1.north east)$) {\pgfplotslegendfromname{simultaneouscombinedlegend}};
\end{tikzpicture}
 \caption{Comparison of the three redistribution methods after two
 simultaneous initial failures with a common surviving neighbor. Curves
 show mean secondary-failure and unserved-load fractions for the grid in
 (a)--(b) and scale-free networks in (c)--(d).}
 \label{fig:simulation-results}
\end{figure}

For the grid network, Fig.~\ref{fig:simulation-results}(a)--(b) shows the results.
Controlled redistribution produces no secondary failures over the tested
range. At the largest tested disturbance, $\Delta=4.5$, it also leaves less
load unserved on average than either of the equal redistribution or the residual capacity redistribution.
Equal redistribution passes on all failed-node load
whenever a surviving neighbor remains, even if that neighbor lacks sufficient
capacity. The resulting overload can cause further failures and if a newly
failed node has no surviving neighbor, its load remains unserved.
Residual-capacity redistribution uses the available spare capacity, but
treats each failed node separately. When several failed nodes share a
receiver, the sum of their transfers may exceed its residual capacity.
Controlled redistribution instead selects all transfers together and
respects each receiver's capacity.

For the scale-free networks, Fig.~\ref{fig:simulation-results}(c)--(d)
shows a tradeoff. Controlled redistribution again produces no secondary
failures. However, at the largest tested disturbance, $\Delta=4.5$, it leaves more
load unserved on average than either heuristic.
Thus, avoiding further failures does not necessarily preserve more served
load. The two optimization problems enforce receiver capacities at the
current step, but do not minimize total unserved load over the complete simulation horizon,
since a redistribution that is optimal at the current step can increase
unserved load later.

\section{Conclusion}
\label{sec:conclusion}

We proposed a topology update approach for irreversible node failures that preserves
row stochasticity and makes unserved load explicit.  We also derived a 
condition for lossless redistribution after simultaneous failures and
formulated a two-stage controller to minimize unserved load and balance
receiver utilization.
For two initial failures with shared receivers, the controller avoided
secondary failures in all tested cases. At the largest disturbance it left
less load unserved than both heuristics on the grid network, but more on the
scale-free networks. The benefit in served load therefore depends on the
network and disturbance, rather than mitigation of failures alone.
Future work will consider
multi-step and robust formulations that account for future cascade stages.

\bibliographystyle{IEEEtran}
\bibliography{references}

@article{motter2002cascade,
  author  = {Adilson E. Motter and Ying-Cheng Lai},
  title   = {Cascade-Based Attacks on Complex Networks},
  journal = {Phys. Rev. E},
  volume  = {66},
  number  = {6},
  pages   = {065102},
  year    = {2002},
  doi     = {10.1103/PhysRevE.66.065102}
}

@article{crucitti2004model,
  author  = {Paolo Crucitti and Vito Latora and Massimo Marchiori},
  title   = {Model for Cascading Failures in Complex Networks},
  journal = {Phys. Rev. E},
  volume  = {69},
  number  = {4},
  pages   = {045104},
  year    = {2004},
  doi     = {10.1103/PhysRevE.69.045104}
}

@article{buldyrev2010catastrophic,
  author  = {Sergey V. Buldyrev and Roni Parshani and Gerald Paul and H. Eugene Stanley and Shlomo Havlin},
  title   = {Catastrophic Cascade of Failures in Interdependent Networks},
  journal = {Nature},
  volume  = {464},
  number  = {7291},
  pages   = {1025--1028},
  year    = {2010},
  doi     = {10.1038/nature08932}
}

@article{barabasi1999emergence,
  author  = {Albert-L\'{a}szl\'{o} Barab\'{a}si and R\'{e}ka Albert},
  title   = {Emergence of Scaling in Random Networks},
  journal = {Science},
  volume  = {286},
  number  = {5439},
  pages   = {509--512},
  year    = {1999},
  doi     = {10.1126/science.286.5439.509}
}

@article{dobson2007complex,
  author  = {Ian Dobson and Benjamin A. Carreras and Vickie E. Lynch and David E. Newman},
  title   = {Complex Systems Analysis of Series of Blackouts: Cascading Failure, Critical Points, and Self-Organization},
  journal = {Chaos},
  volume  = {17},
  number  = {2},
  pages   = {026103},
  year    = {2007},
  doi     = {10.1063/1.2737822}
}

@article{valdez2020cascading,
  author  = {Lucas D. Valdez and Louis Shekhtman and Cristian E. {La Rocca} and Xin Zhang and Sergey V. Buldyrev and Paul A. Trunfio and Lidia A. Braunstein and Shlomo Havlin},
  title   = {Cascading Failures in Complex Networks},
  journal = {J. Complex Netw.},
  volume  = {8},
  number  = {2},
  pages   = {cnaa013},
  year    = {2020},
  doi     = {10.1093/comnet/cnaa013}
}

@inproceedings{bienstock2011optimal,
  author    = {Daniel Bienstock},
  title     = {Optimal Control of Cascading Power Grid Failures},
  booktitle = {Proc. IEEE CDC-ECC},
  address   = {Orlando, FL, USA},
  pages     = {2166--2173},
  year      = {2011},
  doi       = {10.1109/CDC.2011.6160415}
}

@article{como2015throughput,
  author  = {Giacomo Como and Enrico Lovisari and Ketan Savla},
  title   = {Throughput Optimality and Overload Behavior of Dynamical Flow Networks Under Monotone Distributed Routing},
  journal = {IEEE Trans. Control Netw. Syst.},
  volume  = {2},
  number  = {1},
  pages   = {57--67},
  year    = {2015},
  doi     = {10.1109/TCNS.2014.2367361}
}

@article{motter2004control,
  author  = {Adilson E. Motter},
  title   = {Cascade Control and Defense in Complex Networks},
  journal = {Phys. Rev. Lett.},
  volume  = {93},
  number  = {9},
  pages   = {098701},
  year    = {2004},
  doi     = {10.1103/PhysRevLett.93.098701}
}

@article{olfatisaber2007consensus,
  author  = {Reza Olfati-Saber and J. Alex Fax and Richard M. Murray},
  title   = {Consensus and Cooperation in Networked Multi-Agent Systems},
  journal = {Proc. IEEE},
  volume  = {95},
  number  = {1},
  pages   = {215--233},
  year    = {2007},
  doi     = {10.1109/JPROC.2006.887293}
}

@article{kumar2021minimizing,
  author  = {Rajesh Kumar and Suchi Kumari and Manju Bala},
  title   = {Minimizing the Effect of Cascade Failure in Multilayer Networks with Optimal Redistribution of Link Loads},
  journal = {J. Complex Netw.},
  volume  = {9},
  number  = {6},
  pages   = {cnab043},
  year    = {2021},
  doi     = {10.1093/comnet/cnab043}
}

@article{pahwa2013load,
  author  = {Sakshi Pahwa and Caterina Scoglio and Sanjoy Das and Noel Schulz},
  title   = {Load-Shedding Strategies for Preventing Cascading Failures in Power Grid},
  journal = {Electr. Power Compon. Syst.},
  volume  = {41},
  number  = {9},
  pages   = {879--895},
  year    = {2013},
  doi     = {10.1080/15325008.2013.792884}
}

@article{AlTakrouri2013,
  author  = {Saleh Al-Takrouri and Andrey V. Savkin},
  title   = {A Decentralized Flow Redistribution Algorithm for Avoiding Cascaded Failures in Complex Networks},
  journal = {Physica A},
  volume  = {392},
  number  = {23},
  pages   = {6135--6145},
  year    = {2013},
  doi     = {10.1016/j.physa.2013.07.078}
}

@article{Savla2014,
  author  = {Ketan Savla and Giacomo Como and Munther A. Dahleh},
  title   = {Robust Network Routing Under Cascading Failures},
  journal = {IEEE Trans. Netw. Sci. Eng.},
  volume  = {1},
  number  = {1},
  pages   = {53--66},
  year    = {2014},
  doi     = {10.1109/TNSE.2014.2373358}
}

@article{ford1956maximal,
  author  = {Ford, Jr., L. R. and Fulkerson, D. R.},
  title   = {Maximal Flow Through a Network},
  journal = {Can. J. Math.},
  volume  = {8},
  pages   = {399--404},
  year    = {1956},
  doi     = {10.4153/CJM-1956-045-5}
}

\end{document}